\documentclass{amsart}
\usepackage{amssymb}
\usepackage{graphicx}

\newtheorem{thm}{Theorem}
\newtheorem{cor}{Corollary}
\newtheorem{lem}{Lemma}

\newtheorem{rem}{Remark}

\newcommand{\A}{{\mathcal A}}

\newcommand{\U}{{\mathcal U}}
\newcommand{\es}{{\mathcal S}}

\newcommand{\ID}{{\mathbb D}}

\newcommand{\D}{{\mathbb D}}

\def\be{\begin{equation}}
\def\ee{\end{equation}}

\newcommand{\bee}{\begin{enumerate}}
\newcommand{\eee}{\end{enumerate}}

\newcommand{\blem}{\begin{lem}}
\newcommand{\elem}{\end{lem}}
\newcommand{\bthm}{\begin{thm}}
\newcommand{\ethm}{\end{thm}}
\newcommand{\bcor}{\begin{cor}}
\newcommand{\ecor}{\end{cor}}
\newcommand{\beg}{\begin{example}}
\newcommand{\eeg}{\end{example}}
\newcommand{\begs}{\begin{examples}}
\newcommand{\eegs}{\end{examples}}
\newcommand{\bdefe}{\begin{defin}}
\newcommand{\edefe}{\end{defin}}
\newcommand{\bprob}{\begin{prob}}
\newcommand{\eprob}{\end{prob}}
\newcommand{\bei}{\begin{itemize}}
\newcommand{\eei}{\end{itemize}}

\newcommand{\bcon}{\begin{conj}}
\newcommand{\econ}{\end{conj}}
\newcommand{\bcons}{\begin{conjs}}
\newcommand{\econs}{\end{conjs}}
\newcommand{\bprop}{\begin{propo}}
\newcommand{\eprop}{\end{propo}}
\newcommand{\br}{\begin{rem}}
\newcommand{\er}{\end{rem}}
\newcommand{\brs}{\begin{rems}}
\newcommand{\ers}{\end{rems}}
\newcommand{\bo}{\begin{obser}}
\newcommand{\eo}{\end{obser}}
\newcommand{\bos}{\begin{obsers}}
\newcommand{\eos}{\end{obsers}}
\newcommand{\bpf}{\begin{pf}}
\newcommand{\epf}{\end{pf}}
\newcommand{\ba}{\begin{array}}
\newcommand{\ea}{\end{array}}
\newcommand{\beq}{\begin{eqnarray}}
\newcommand{\beqq}{\begin{eqnarray*}}
\newcommand{\eeq}{\end{eqnarray}}
\newcommand{\eeqq}{\end{eqnarray*}}

\begin{document}
\bibliographystyle{amsplain}

\title[Hankel determinants of second and third order]{Hankel determinants of second and third order for the class $\boldsymbol{\mathcal{S}}$ of univalent functions}

\author[M. Obradovi\'{c}]{Milutin Obradovi\'{c}}
\address{Department of Mathematics,
Faculty of Civil Engineering, University of Belgrade,
Bulevar Kralja Aleksandra 73, 11000, Belgrade, Serbia}
\email{obrad@grf.bg.ac.rs}

\author[N. Tuneski]{Nikola Tuneski}
\address{Department of Mathematics and Informatics, Faculty of Mechanical Engineering, Ss. Cyril and Methodius
University in Skopje, Karpo\v{s} II b.b., 1000 Skopje, Republic of North Macedonia.}
\email{nikola.tuneski@mf.edu.mk}







\begin{abstract}
In this paper we give the upper bounds of the Hankel determinants of the second and third order for the class  $\mathcal{S}$ of univalent functions in the unit disc.
\end{abstract}

\maketitle

Let $\mathcal{A}$ be the class of functions $f$ that are analytic  in the open unit disc $\D=\{z:|z|<1\}$ of the form $f(z)=z+a_2z^2+a_3z^3+\cdots$ and let $\mathcal{S}$ be the class of univalent functions in the unit disc $\D$.
Let $\mathcal{S}^{\star}$ and $\mathcal{K}$ denote the subclasses of
${\mathcal A}$ which are starlike and convex in $\ID$, respectively, and
let $\mathcal{U} $ denote the set of all  $f\in {\mathcal A}$ in $\ID$ satisfying
the condition

$$\left |\left (\frac{z}{f(z)} \right )^{2}f'(z)-1\right | < 1 \quad\quad (z\in \ID).$$
(see  \cite{OP-01}, \cite{OP_2011}, \cite{OP_2016}).

\medskip

The $qth$ Hankel determinant for a function $f$ from $\A$ is defined for $q\geq 1$, and
$n\geq 1$ by
\[
        H_{q}(n) = \left |
        \begin{array}{cccc}
        a_{n} & a_{n+1}& \ldots& a_{n+q-1}\\
        a_{n+1}&a_{n+2}& \ldots& a_{n+q}\\
        \vdots&\vdots&~&\vdots \\
        a_{n+q-1}& a_{n+q}&\ldots&a_{n+2q-2}\\
        \end{array}
        \right |.
\]
Thus, the second Hankel determinant is 
\begin{equation}\label{eq_1}
H_{2}(2)= a_2a_4-a_{3}^2
\end{equation}
and the third is
\[ H_3(1) =  \left |
        \begin{array}{ccc}
        1 & a_2& a_3\\
        a_2 & a_3& a_4\\
        a_3 & a_4& a_5\\
        \end{array}
        \right | = a_3(a_2a_4-a_{3}^2)-a_4(a_4-a_2a_3)+a_5(a_3-a_2^2).
\]
The concept of Hankel determinant finds its application in the theory of singularities (see \cite{dienes}) and in the study of power series with integral coefficients.

\medskip

For some subclasses of the class $\mathcal{S}$ of univalent functions the sharp estimation of $|H_{2}(2)|$ are known. For example,
for the classes $\mathcal{S}^{\star}$ and $\mathcal{U} $ we have that $|H_{2}(2)|\leq 1$
(see \cite{JHD,OT_2019}), while $|H_{2}(2)|\leq \frac{1}{8} $ for the class $\mathcal{K}$ (\cite{JHD}).
Finding sharp estimates of the third order Hankel determinant turns out to be more complicated, so very few are known.
An overview of results on the upper bound of $|H_3(1)|$ can be found in \cite{shi}, while new non-sharp upper bounds for different classes and conjectures about the sharp ones are given in \cite{OT_2019-3}.

\medskip

In this paper we give an upper bound of $|H_{2}(2)|$ and $|H_{3}(1)|$
for the class $\mathcal{S}$. Namely, we have

\bthm\label{19-th 1} For the class $\mathcal{S}$ we have
\[
|H_{2}(2)|\leq A , \quad\mbox{where}\quad 1\leq A\leq \frac{11}{3}=3,66\ldots
\]
and
\[
|H_{3}(1)|\leq  B, \quad\mbox{where}\quad \frac49\leq B\leq \frac{32+\sqrt{285}}{15} = 3.258796\cdots
\]
\ethm

\medskip

\begin{proof}
In the proof of this theorem we will use mainly the notations and results given in the book of N. A. Lebedev (\cite{Lebedev}).

Let $f \in \mathcal{S}$ and let
\[
\log\frac{f(t)-f(z)}{t-z}=\sum_{p,q=0}^{\infty}\omega_{p,q}t^{p}z^{q},
\]
where $\omega_{p,q}$ are called Grunsky's coefficients with property $\omega_{p,q}=\omega_{q,p}$.
For those coefficients we have the next Grunsky's inequality (\cite{duren,Lebedev}):
\be\label{eq 4}
\sum_{q=1}^{\infty}q \left|\sum_{p=1}^{\infty}\omega_{p,q}x_{p}\right|^{2}\leq \sum_{p=1}^{\infty}\frac{|x_{p}|^{2}}{p},
\ee
where $x_{p}$ are arbitrary complex numbers such that last series converges.

Further, it is well-known that if
\be\label{eq 5}
f(z)=z+a_{2}z^{2}+a_{3}z^{3}+...
\ee
belongs to $\mathcal{S}$, then also
\[
f_{2}(z)=\sqrt{f(z^{2})}=z +c_{3}+c_{5}z^{5}+...
\]
belongs to the class $\mathcal{S}$. Then for the function $f_{2}$ we have the appropriate Grunsky's
coefficients of the form $\omega_{2p-1,2q-1}^{(2)}$ and the inequality \eqref{eq 4} has the form
\be\label{eq 7}
\sum_{q=1}^{\infty}(2q-1) \left|\sum_{p=1}^{\infty}\omega_{2p-1,2q-1}^{(2)}x_{2p-1}\right|^{2}\leq \sum_{p=1}^{\infty}\frac{|x_{2p-1}|^{2}}{2p-1}.
\ee
As it has been shown in \cite[p.57]{Lebedev}, if $f$ is given by \eqref{eq 5} then the coefficients $a_{2}$, $a_{3}$, $a_{4}$ and $a_5$
are expressed by Grunsky's coefficients  $\omega_{2p-1,2q-1}^{(2)}$ of the function $f_{2}$ given by
\eqref{eq 5} in the following way (in the next text we omit upper index 2 in $\omega_{2p-1,2q-1}^{(2)}$):

\be\label{eq 8}
\begin{split}
a_{2}&=2\omega _{11},\\
a_{3}&=2\omega_{13}+3\omega_{11}^{2}, \\
a_{4}&=2\omega_{33}+8\omega_{11}\omega_{13}+\frac{10}{3}\omega_{11}^{3}\\
a_{5}&=2\omega_{35}+8\omega_{11}\omega_{33}+5\omega_{15}^{2}+18\omega_{11}^{2}\omega_{13}+\frac{7}{3}\omega_{11}^{4}\\
0&=3\omega_{15}-3\omega_{11}\omega_{13}+\omega_{11}^{3}-3\omega_{33}.
\end{split}
\ee
Now, from \eqref{eq_1} and \eqref{eq 8} we have
\[
\begin{split}
H_{2}(2)&= 4\omega_{11}\omega_{33}+4\omega_{11}^{2}\omega_{13}-4\omega_{13}^{2}-\frac{7}{3}\omega_{11}^{4}\\
&= 4\omega_{11}\omega_{33}-\frac{4}{3}\omega_{11}^{4}-\left (2\omega_{13}-\omega_{11}^{2}\right)^{2},
\end{split}
\]
and from here
\be\label{eq 9}
|H_{2}(2)|\leq 4|\omega_{11}||\omega_{33}|+\frac{4}{3}|\omega_{11}|^{4}+\left|2\omega_{13}-\omega_{11}^{2}\right|^{2}.
\ee
Since for the class $\mathcal{S}$ we have $ |a_{3}-a_{2}^{2}|\leq1$ (see \cite{duren}) and since from \eqref{eq 8}
$$ |2\omega_{13}-\omega_{11}^{2}| =|a_{3}-a_{2}^{2}|,$$
then
\be\label{eq 10}
|2\omega_{13}-\omega_{11}^{2}|\leq1 .
\ee
On the other hand, from \eqref{eq 7} for $x_{2p-1}=0$, $p=3,4,\ldots$ we have
\be\label{eq 11}
|\omega_{11}x_{1}+\omega_{31}x_{3}|^{2}+3|\omega_{13}x_{1}+\omega_{33}x_{3}|^{2}
\leq |x_{1}|^{2}+\frac{|x_{3}|^{2}}{3}.
\ee
From \eqref{eq 11} for $x_{1}=1,\,x_{3}=0 $ and since $\omega_{31}=\omega_{13}$, we have
$$|\omega_{11}|^{2}+3|\omega_{13}|^{2}\leq1 ,$$
which implies
\be\label{eq 12}
|\omega_{13}|^{2}\leq\frac{1}{3}(1-|\omega_{11}|^{2}).
\ee
Also, for $x_{1}=0,\,x_{3}=1 $ we obtain
$$|\omega_{31}|^{2}+3|\omega_{33}|^{2}\leq\frac{1}{3} $$
and so
\be\label{eq 13}
|\omega_{33}|\leq\frac13\sqrt{1-3|\omega_{31}|^2}\le \frac{1}{3}.
\ee

Finally, from \eqref{eq 9}, \eqref{eq 10}, \eqref{eq 12} and \eqref{eq 13} we have
\[
|H_{2}(2)|\leq \frac{4}{3}|\omega_{11}|+\frac{4}{3}|\omega_{11}|^{4}+1  \leq \frac{11}{3},
\]
because from \eqref{eq 8} we have that
$$|a_{2}|=|2\omega_{11}|\leq2 \quad\Rightarrow \quad|\omega_{11}|\leq 1 .$$

Since $\es^\ast$ and $\U$ are both subsets of $\es$ with 1 as a sharp upper bound of $|H_{2}(2)|$, we have that on the class $\es$, $|H_{2}(2)|\ge1$.

\medskip

As for Hankel determinant of the third order, by using \eqref{eq 8}, we can write
\[
\begin{split}
H_3(1) &= a_3(a_2a_4-a_{3}^2)-a_4(a_4-a_2a_3)+a_5(a_3-a_2^2)\\
&= -8\omega_{13}^{3}+2\omega_{11}^{4}\omega_{13}+\frac{8}{3}\omega_{11}^{3}\omega_{33}-4\omega_{33}^{2} -\frac{4}{9}\omega_{11}^{6} \\
& \quad +4\omega_{13}\omega_{35}+10\omega_{13}\omega_{15}^{2}-5\omega_{11}^{2}\omega_{15}^{2} -2\omega_{11}^2\omega_{35} \\
&=-2\omega_{13}\left(4\omega_{13}^{2}-\omega_{11}^{4}\right)-\left(2\omega_{33}-\frac{2}{3}\omega_{11}^{3}\right)^{2}+
(2\omega_{35}+5\omega_{15}^{2})(2\omega_{13}-\omega_{11}^{2}),
\end{split}
\]
and from here
\[
\begin{split}
|H_{3}(1)|&\le \underbrace{2|\omega_{13}|\left|4\omega_{13}^{2}-\omega_{11}^{4}\right|}_{B_1}+\underbrace{\left|2\omega_{33}-\frac{2}{3}\omega_{11}^{3}\right|^{2}}_{B_2}+
\underbrace{|2\omega_{35}+5\omega_{15}^{2}||(2\omega_{13}-\omega_{11}^{2}|}_{B_3}\\
&= B_{1}+B_{2}+B_{3}.
\end{split}
\]

By using the relations \eqref{eq 10} and \eqref{eq 12}, we obtain
\[
\begin{split}
B_{1}&= 2\cdot |\omega_{13}|\cdot \left|2\omega_{13}-\omega_{11}^{2}\right|\cdot \left|2\omega_{13}+\omega_{11}^{2}\right|\\
&\leq 2\cdot |\omega_{13}|\cdot \left|2\omega_{13}+\omega_{11}^{2}\right|\\
&\leq 2\cdot |\omega_{13}|\cdot \left(2|\omega_{13}|+|\omega_{11}|^{2}\right)\\
& = 4|\omega_{13}|^2+2|\omega_{13}| \cdot |\omega_{11}|^2\\
&\leq \frac{2}{3}\cdot \left[2\left(1-|\omega_{11}|^{2}\right)+\sqrt{3}\cdot |\omega_{11}|^{2}\cdot \sqrt{1-|\omega_{11}|^{2}}\right]\\
&=:\frac{2}{3}\cdot \varphi(|\omega_{11}|^{2}),
\end{split}
\]
where
$$ \varphi(t)=2(1-t)+\sqrt{3}t\sqrt{1-t},\quad\quad0\leq t\leq1.$$

It is easily to show that the function $\varphi$ decreases on $(0,1)$ and has maximum $\varphi(0)=2$, which implies
\be\label{eq 14}
B_{1}\leq \frac{2}{3}\cdot \varphi(0)=\frac43.
\ee

From the last equation in the relation \eqref{eq 8}, we have
$$2\omega_{33}-\frac{2}{3}\omega_{11}^{3}=2\omega_{15}-2\omega_{11}\omega_{13},$$
and from here
\be\label{eq 15}
\left|2\omega_{33}-\frac{2}{3}\omega_{11}^{3}\right|\leq2|\omega_{15}|+2|\omega_{11}||\omega_{13}|.
\ee

Similarly as in \eqref{eq 11}, we have
\be\label{eq 16}
|\omega_{11}x_{1}+\omega_{31}x_{3}|^{2}+3|\omega_{13}x_{1}+\omega_{33}x_{3}|^{2}+
5|\omega_{15}x_{1}+\omega_{35}x_{3}|^{2}
\leq |x_{1}|^{2}+\frac{|x_{3}|^{2}}{3}.
\ee
If we put $x_{1}=1$ and $x_{3}=0$, then we get
$$|\omega_{11}|^{2}+3|\omega_{13}|^{2}+5|\omega_{15}|^{2}\leq 1,$$
and so
\be\label{eq 17}
|\omega_{15}|\leq \frac{1}{\sqrt{5}}\sqrt{1-|\omega_{11}|^{2}-3|\omega_{13}|^{2}}.
\ee
From \eqref{eq 15} and \eqref{eq 17}, we have
\[
\begin{split}
\left|2\omega_{33}-\frac{2}{3}\omega_{11}^{3}\right| &\leq \frac{2}{\sqrt{5}} \left(\sqrt{1-|\omega_{11}|^{2}-3|\omega_{13}|^{2}}+\sqrt{5} |\omega_{11}|\cdot |\omega_{13}|\right)\\
&=: \frac{2}{\sqrt{5}}\cdot \psi(|\omega_{11}|,|\omega_{13}|),
\end{split}
\]
where
$$\psi(t,s)=\sqrt{1-t^{2}-3s^{2}}+\sqrt{5}ts,\quad\quad 0\leq t\leq1, \quad 0\leq s \leq \frac{1}{\sqrt{3}}\sqrt{1-t^{2}}.$$
It is an elementary fact to find that in cited domain $\max \psi =1$ attained for $t=s=0$, which implies
\be\label{eq 18}
B_{2}=\left|2\omega_{33}-\frac{2}{3}\omega_{11}^{3}\right|^{2}\leq \left( \frac{2}{\sqrt{5}}\right)^{2}=\frac{4}{5}.
\ee

From relation \eqref{eq 16} we also have
$$5|\omega_{15}x_{1}+\omega_{35}x_{3}|^{2} \leq |x_{1}|^{2}+\frac{|x_{3}|^{2}}{3}.$$
If we put in the previous relation $x_{1}=5\omega_{15}$, $x_{3}=2$, and then use \eqref{eq 17} we receive
$$|2\omega_{35}+5\omega_{15}^{2}|^{2}\leq 5|\omega_{15}|^{2}+\frac{4}{15}\leq 1-|\omega_{11}|^{2}-3|\omega_{13}|^{2}+\frac{4}{15}\leq \frac{19}{15},$$
which finally gives
\be\label{eq 19}
B_{3}=|2\omega_{35}+5\omega_{15}^{2}|\cdot |2\omega_{13}-\omega_{1}^{2}|\leq \sqrt{\frac{19}{15}}
\ee
(in the last step we have used the relation \eqref{eq 10}). By using relations \eqref{eq 14}, \eqref{eq 18} and \eqref{eq 19}, we obtained
$$|H_{3}(1)|\leq B_{1}+B_{2}+B_{3}\leq \frac43+\frac45+\sqrt{\frac{19}{15}} = \frac{32+\sqrt{285}}{15}.$$

The function defined by $\frac{zf'(z)}{f(z)}=\frac{1+z^3}{1-z^3}$ where $a_2=a_3=a_5=0$, $a_4=\frac23$ is starlike (thus univalent) and $H_3(1)=-\frac49$. Therefore on the class $\es$,
$$|H_{3}(1)|\ge \frac49.$$
\end{proof}

\end{document}